\documentclass[10pt,twoside,reqno]{amsart}
\usepackage{amssymb}
\usepackage{multirow}
\usepackage{graphicx}
\usepackage[dvips]{epsfig}
\usepackage{latexsym}
\usepackage{amsmath}
\usepackage{amsthm}
\usepackage{amscd}
\usepackage{amssymb}
\input xy
\xyoption{all}

\calclayout

\numberwithin{equation}{section}
\makeatletter

\renewcommand{\@secnumfont}{\bfseries}

\renewcommand{\section}{\@startsection{section}{1}%
  {0mm}{.7\linespacing\@plus\linespacing}{.5\linespacing}
  {\normalfont\bfseries\centering}}

\newcommand{\bibsection}{\@startsection{section}{1}%
  {0mm}{.7\linespacing\@plus\linespacing}{.5\linespacing}
  {\normalfont\scshape\centering}}

\renewcommand{\@biblabel}[1]{#1.}

\newtheorem{thm}{\bf Theorem}[section]
\newtheorem{lem}[thm]{\bf Lemma}

\newtheorem{cor}[thm]{\bf Corollary}

\begin{document}

\vspace{1.3cm}

\title {On degenerate gamma functions}

\author{Taekyun  Kim${}^{1}$}
\address{${}^{1}$Department of Mathematics, Kwangwoon University, Seoul 139-701, Republic of Korea}
\email{tkkim@kw.ac.kr}
\author{Dae San  Kim${}^{2}$}
\address{${}^{2}$Department of Mathematics, Sogang University, Seoul 121-742, Republic of Korea} \email{dskim@sogang.ac.kr}


\maketitle

\begin{abstract} Recently, the degenerate gamma functions are introduced as a degenerate version of the usual gamma function by Kim--Kim. In this paper, we investigate several properties of them. Namely, we obtain an analytic continuation as a meromorphic function on the whole complex plane, the difference formula, the values at positive integers, some expressions following from the Weierstrass and Euler formulas for the usual gamma function and an integral representation as the integral along a Hankel contour.
\end{abstract}

\pagestyle{myheadings}
\markboth{\centerline {\scriptsize T. Kim, D. S. Kim}}
         {\centerline{\scriptsize On degenerate gamma functions}}
\bigskip
\medskip
\section{\bf Introduction}
\medskip

It is not an exaggeration to say that the gamma function is the most important nonelementary transcendental function. It appears in many areas such as hypergeometric series, asymptotic series, definite integration, Riemann zeta function, $L$-functions and number theory in general.
The gamma function was introduced by Euler and subsequently studied by eminent
mathematicians like Daniel Bernoulli, Legendre, Gauss, Liouville, Weierstrass,
Hermite as well as many other people.\\
\indent In [6], the degenerate gamma functions were introduced in an attempt to find a
degenerate version of the usual gamma function and the related degenerate
Laplace transform was also studied (see also [7]).
In this paper, we would like to derive some basic properties about the degenerate
gamma functions including analytic continuation to a meromorphic function to the whole complex plane, the difference formula, some expressions comimg from the Euler and Weierstrass formulas, integral representation as an integal along a Hankel contour. Also, we will discuss the degenerate Beta functions.\\
\indent In the rest of this section, we will recall some basic facts on the gamma function, the beta function, the degenerate exponential function and the degenerate gamma functions. In Section 2, we will study our
main results of this paper. Finally, we will discuss the integral representation of the degenerate gamma functions in Section 3 and conclude our results in Section 4.

\vspace{1mm}

For $s \in \mathbb{C}$ with $\textnormal{Re}(s) > 0$, the gamma function is defined by\\
\begin{equation}\label{01}
\begin{split}
\Gamma(s) = \int_{0}^{\infty} e^{-t} t^{s-1} dt,\,\,\,\, (\textnormal{see}\,\,[1-22]).
\end{split}
\end{equation}
From \eqref{01}, we note that\\
\begin{equation}\label{01-1}
\begin{split}
\Gamma(s+1) = s\Gamma(s),\,\,\, \Gamma(n) = (n-1)!, \,\,\,\, (n \in \mathbb{N}), \,\,\,\,(\textnormal{see}\,\, [15,16,17,18]).
\end{split}
\end{equation}
For $\alpha, \beta \in \mathbb{C}$ with $\textnormal{Re}(\alpha)> 0$ and $\textnormal{Re}(\beta)> 0$, the Beta function is defined by the integrals:
\begin{equation}\label{02}
\begin{split}
B(\alpha, \beta) = \int^{1}_{0} t^{\alpha -1}(1-t)^{\beta -1} dt= \int_{0}^{\infty}\frac{u^{\alpha-1}}{(1+u)^{\alpha+\beta}} du, \,\,\,\, (\textnormal{see}\,\,[1,15,16,17]).
\end{split}
\end{equation}
From \eqref{01} and \eqref{02}, we note that
\begin{equation}\label{03}
B(\alpha ,\beta)=\frac{\Gamma(\alpha) \Gamma(\beta)}{\Gamma(\alpha + \beta)},  \,\,\,\,(\textnormal{see}\,\,[15,22]).
\end{equation}

\vspace{1mm}

As is well known, the Euler's formula for the gamma function is given by
\begin{equation}\label{04}
\begin{split}
\Gamma(z)  & = \frac{1}{z} \prod^{\infty}_{n=1} \left \{\left(1+\frac{1}{n}\right)^{z}\left(1+\frac{z}{n}\right)^{-1} \right\}\\
& = \lim_{n \to \infty} \frac{(n-1)!}{z(z+1) \cdots (z+n-1)}n^{z} ,\,\,\,\, (\textnormal{see}\,\,[2,3,4,13,16,22]),
\end{split}
\end{equation}
where $z \neq 0,-1,-2,-3,\cdots .$
From \eqref{04}, we can derive the Weierstrass formula which is given by \\
\begin{equation}\label{05}
\begin{split}
\frac{1}{\Gamma(z)} =ze^{\gamma z} \prod^{\infty}_{n=1} \left \{ \left(1+\frac{z}{n}\right)e^{-\frac{z}{n}}\right \}, \,\,\,\,(\textnormal{see}\,\,[2,3,4,22]).
\end{split}
\end{equation}
where $\gamma$ is Euler constant defined by\\
\begin{equation*}
\begin{split}
\gamma = \lim_{n \to \infty} \left(1+\frac{1}{2}+\frac{1}{3} \cdots + \frac{1}{n}-\log n\right).
\end{split}
\end{equation*}
\indent We note here that we may speak of the analytically continued beta function obtained from \eqref{03} which is analytic except for $\alpha=0,-1,-2,\dots, \beta=0,-1,-2,\dots$. Indeed, from \eqref{03} and \eqref{05}, we have
\begin{equation*}
B(\alpha,\beta)= \frac{\Gamma(\alpha)\Gamma(\beta)}{\Gamma(\alpha+\beta)}=\frac{\alpha+\beta}{\alpha\beta}\prod_{n=1}^{\infty}\frac{(1+\frac{\alpha+\beta}{n})}{(1+\frac{\alpha}{n})(1+\frac{\beta}{n})}.
\end{equation*}
From \eqref{05}, we note that
\begin{equation}\label{06}
\begin{split}
\Gamma(z)\Gamma(-z) &= -\frac{1}{z^{2}} \prod^{\infty}_{n=1} \left\{\left(1+\frac{z}{n}\right) e^{-\frac{z}{n}}\right\}^{-1}\prod^{\infty}_{n=1}\left\{ \left(1-\frac{z}{n}\right)e^{\frac{z}{n}}\right\}^{-1} \\
&= \frac{-\pi}{z \sin \pi z}.
\end{split}
\end{equation}\\
Thus, by \eqref{06} and \eqref{01-1}, we get
\begin{equation}\label{07}
\begin{split}
\Gamma(z)\Gamma(1-z) = \Gamma(z)(-z\Gamma(-z)) = \frac{\pi}{\sin \pi z} ,\,\,\,\,(\textnormal{see}\,\,[16]).
\end{split}
\end{equation}

\vspace{1mm}
For any nonzero $\lambda \in \mathbb{R}$, the degenerate exponential functions are defined as
\begin{equation}\label{08}
\begin{split}
e_{\lambda}^{x}(t) = (1+\lambda t)^{\frac{x}{\lambda}},\,\,\,  e_{\lambda}(t)=e_{\lambda}^{1}(t) =(1+\lambda t)^{\frac{1}{\lambda}} , \,\,\,\,(\textnormal{see}\,\,[6,7,8,9,10]).
\end{split}
\end{equation}
Note that $\lim_{\lambda \to 0} e^{x}_{\lambda}(t) = e^{xt}$.\\
Let $\log_{\lambda} (t)$ be the compositional inverse of $e_{\lambda}(t)$, which is called the degenerate logarithm. Then we have \\
\begin{equation}\label{09}
\begin{split}
\log_{\lambda}(t) = \frac{1}{\lambda}\left(t^{\lambda}-1\right), \,\,\,\,(\textnormal{see}\,\,[8,9,10]).
\end{split}
\end{equation}
From \eqref{08} and \eqref{09}, we note that\\
\begin{equation*}
\begin{split}
t = \log_{\lambda}e_{\lambda}(t) = e_{\lambda}\log_{\lambda}(t).
\end{split}
\end{equation*}
For $\lambda \in (0,1)$, the degenerate gamma function for complex variable $s$ with $0<\textnormal{Re}(s)<\frac{1}{\lambda}$
is defined by Kim-Kim as\\
\begin{equation}\label{10}
\begin{split}
\Gamma_{\lambda}(s) = \int^{\infty}_{0} e^{-1}_{\lambda}(t)t^{s-1} dt = \int^{\infty}_{0}(1+\lambda t)^{-\frac{1}{\lambda}}t^{s-1}dt,\,\,\,\,(\textnormal{see}\,\,[6]).
\end{split}
\end{equation}
Thus, by \eqref{10}, we get\\
\begin{equation}\label{11}
\begin{split}
\Gamma_{\lambda}(s+1) = \frac{s}{(1-\lambda)^{s+1}} \Gamma_{\frac{\lambda}{1-\lambda}}(s)\,\,\,\, 0<\textnormal{Re}(s)<\frac{1-\lambda}{\lambda}, \,\,\,\,(\textnormal{see}\,\,[6,19,20]).
\end{split}
\end{equation}
From\eqref{11}, we note that, for any integer $k \geq 0$, we have\\
\begin{equation}\label{12}
\begin{split}
\Gamma_{\lambda}(s+1) &= \frac{s(s-1)(s-2)\cdots (s-(k+1)+1)}{(1-\lambda)(1-2\lambda)\cdots (1-k\lambda)(1-(k+1)\lambda)^{s-k+1}}\\
&\quad    \times \Gamma_{\frac{\lambda}{1-(k+1)\lambda}}(s-k), \,\,\,\,(\textnormal{see}\,\,[6]),
\end{split}
\end{equation}
where $\lambda \in ( 0, \frac{1}{k+1})$ with $k \in \mathbb{N}$ and $k < \textnormal{Re}(s)<\frac{1-\lambda}{\lambda}$.
In particular, for $k \in \mathbb{N}$ and $\lambda \in (0,\frac{1}{k})$, we have\\
\begin{equation}\label{13}
\begin{split}
\Gamma_{\lambda}(k) = \frac{(k-1)!}{(1-\lambda)(1-2\lambda)\cdots (1-k\lambda)}, \,\,\,\,(\textnormal{see}\,\,[6]).
\end{split}
\end{equation}\\

\medskip

\section{\bf Degenerate gamma functions}
\medskip

In this section, we assume that $\lambda$ is any fixed real number with $0 < \lambda <1$.
By \eqref{02}, we have\\
\begin{equation}\label{16}
\begin{split}
\Gamma_{\lambda}(s) =  \lambda^{-s}\int^{\infty}_{0}\frac{t^{s-1}}{(1+t)^{\frac{1}{\lambda}}}dt = \lambda^{-s} B\left(s, \frac{1}{\lambda}-s\right) = \lambda^{-s}\frac{\Gamma(s)\Gamma(\frac{1}{\lambda}-s)}{\Gamma(\frac{1}{\lambda})}.
\end{split}
\end{equation}
We note here that \eqref{16} holds initially for $s \in \mathbb{C}$ with $0<\textnormal{Re}(s)<\frac{1}{\lambda}$ and that it further holds for any $s \in \mathbb{C}-\left\{ 0,-1,-2, \dots, \frac{1}{\lambda}, \frac{1}{\lambda}+1,\frac{1}{\lambda}+2, \dots \right\}$ by analytic continuation where it defines an analytic function. Moreover, it has simple poles at $s=-n, \,(n=0,1,2.\dots)$ with residues \\
\begin{equation}\label{16-1}
\textnormal{res}_{s=-n}=\frac{(-1)^n\lambda^n \Gamma(\frac{1}{\lambda}+n)}{n!\Gamma(\frac{1}{\lambda})},
\end{equation}
and simple poles at $s=\frac{1}{\lambda}+n, \,(n=0,1,2.\dots)$ with residues \\
\begin{equation}\label{16-2 }
\textnormal{res}_{s=\frac{1}{\lambda}+n}=\frac{(-1)^{n-1}\lambda^{-n-\frac{1}{\lambda}} \Gamma(\frac{1}{\lambda}+n)}{n!\Gamma(\frac{1}{\lambda})}.
\end{equation}
Therefore, by \eqref{16}, we obtain the following lemma.
\begin{lem}\label{Lemma 1}
	Let $\lambda \in (0,1)$. Then $\Gamma_{\lambda}(s)$ initially defined for  $0<\textnormal{Re}(s)<\frac{1}{\lambda}$ has an analytic continuation to a meromorphic function on $\mathbb{C}$ whose only singularities are simple poles at $s= 0,-1,-2, \dots, \frac{1}{\lambda}, \frac{1}{\lambda}+1,\frac{1}{\lambda}+2, \dots$, with residues \\
\begin{equation*}
\textnormal{res}_{s=-n}=\frac{(-1)^n\lambda^n \Gamma(\frac{1}{\lambda}+n)}{n!\Gamma(\frac{1}{\lambda})}, \,\,
\textnormal{res}_{s=\frac{1}{\lambda}+n}=\frac{(-1)^{n-1}\lambda^{-n-\frac{1}{\lambda}} \Gamma(\frac{1}{\lambda}+n)}{n!\Gamma(\frac{1}{\lambda})}.
\end{equation*}
In addition, for any $s \in \mathbb{C}-\left\{0,-1,-2, \dots, \frac{1}{\lambda}, \frac{1}{\lambda}+1,\frac{1}{\lambda}+2, \dots \right\}$ the following identities hold:
\begin{equation*}
\Gamma_{\lambda}(s) = \lambda^{-s}\frac{\Gamma(s)\Gamma(\frac{1}{\lambda}-s)}{\Gamma(\frac{1}{\lambda})}=\lambda^{-s}B(s,\frac{1}{\lambda}-s),
\end{equation*}
\begin{equation*}
\lambda^s\Gamma_{\lambda}(s)=\lambda^{\frac{1}{\lambda}-s}\Gamma_{\lambda}(\frac{1}{\lambda}-s).
\end{equation*}
\end{lem}

From  Lemma 2.1, we note that\\
\begin{equation}\label{17}
\begin{split}
\Gamma_{\lambda}(s+1) & = \lambda^{-s-1}\frac{\Gamma(s+1)\Gamma(\frac{1}{\lambda}-s-1)}{\Gamma(\frac{1}{\lambda})} \\
&= \frac{s\lambda^{-s}}{\lambda(\frac{1}{\lambda}-s-1)}\frac{\Gamma(s)\left(\frac{1}{\lambda}-s-1\right)\Gamma\left(\frac{1}{\lambda}-s\right)}{\Gamma\left(\frac{1}{\lambda}\right)}\\
& = \frac{s}{\lambda \left(\frac{1}{\lambda}-s-1\right)}\lambda^{-s}\frac{\Gamma(s)\Gamma(\frac{1}{\lambda}-s)}{\Gamma\left(\frac{1}{\lambda}\right)} = \frac{s}{1-\lambda(s+1)}\Gamma_{\lambda}(s).
\end{split}
\end{equation}

From \eqref{10} and \eqref{17}, we obtain the following theorem.\\
\begin{thm}\label{Theorem 2}
	Let $\lambda \in (0,1)$. Then the following difference equation holds for any  $s \in \mathbb{C}-\left\{0,-1,-2\dots, \frac{1}{\lambda}-1, \frac{1}{\lambda},\frac{1}{\lambda}+1, \dots \right\}$\textnormal{:}
	\begin{equation*}
	\begin{split}
	\Gamma_{\lambda}(s+1) = \frac{s}{1-\lambda(s+1)}\Gamma_{\lambda}(s).
	\end{split}
	\end{equation*}
\end{thm}

Let $s=k$ be a positive integer. Then, by Theorem 2.2, we get
\begin{equation}\label{18}
\begin{split}
\Gamma_{\lambda}(k+1)& = \frac{k}{1-\lambda(k+1)}\Gamma_{\lambda}(k) = \frac{k(k-1)}{(1-\lambda(k+1))(1-\lambda k)}\Gamma_{\lambda}(k-1)\\
&=\cdots  = \frac{k\cdot(k-1)\cdots 2\cdot1}{(1-(k+1)\lambda)(1-k\lambda)\cdots(1-2\lambda)}\Gamma_{\lambda}(1).
\end{split}
\end{equation}

From \eqref{16}, we have\\
\begin{equation}\label{19}
\Gamma_{\lambda}(1) = \frac{1}{\lambda}\frac{\Gamma(1)\Gamma\left(\frac{1}{\lambda}-1\right)}{\Gamma\left(\frac{1}{\lambda}\right)} = \frac{1}{\lambda\left(\frac{1}{\lambda}-1\right)}=\frac{1}{1-\lambda}.
\end{equation}

Thus, by \eqref{18} and \eqref{19}, for any integer $k \geq 0$ we get
\begin{equation}\label{20}
\Gamma_{\lambda}(k+1) = \frac{\Gamma(k+1)}{(1)_{k+2,\lambda}},
\end{equation}
where $(x)_{0,\lambda} = 1$, $(x)_{n,\lambda} = x(x-\lambda)\cdots (x-(n-1)\lambda)$, $(n\ge 1)$, (\textnormal{see}\,\,[5-11]). \\

Therefore, by \eqref{20}, we obtain the following theorem.
\begin{thm}{\label{Theorem 3}}
Let $\lambda \in (0,1)$. Then, for $k \in \mathbb{N}$, we have
\begin{equation*}
\begin{split}
\Gamma_{\lambda}(k) = \frac{\Gamma(k)}{(1)_{k+1,\lambda}}.
\end{split}
\end{equation*}
\end{thm}

From \eqref{05} and Lemma 2.1, we note that
\begin{equation}\label{21}
\begin{split}
\frac{1}{\Gamma_{\lambda}(z)} & = \Gamma\left(\frac{1}{\lambda}\right)\lambda^{z}\frac{1}{\Gamma(z)}\frac{1}{\Gamma\left(\frac{1}{\lambda}-z\right)}\\
& = \lambda^{z}\Gamma\left(\frac{1}{\lambda}\right)ze^{\gamma z}\prod^{\infty}_{n=1} \left\{\left(1+\frac{z}{n}\right)e^{-\frac{z}{n}}\right\}\left(\frac{1}{\lambda}-z\right)e^{\gamma\left(\frac{1}{\lambda}-z\right)}\\
&\quad \times \prod^{\infty}_{n=1}\left\{\left(1+\frac{\frac{1}{\lambda}-z}{n}\right)e^{-\frac{\frac{1}{\lambda}-z}{n}} \right\}\\
& = \lambda^{z}\Gamma\left(\frac{1}{\lambda}\right)z\left(\frac{1}{\lambda}-z\right)e^{\frac{\gamma}{\lambda}}\prod^{\infty}_{n=1}\left\{ e^{-\frac{1}{n\lambda}}\left(1+\frac{z}{n}\right)\left( 1+\frac{\frac{1}{\lambda}-z}{n}\right)\right\}.
\end{split}
\end{equation}
Also, from this equation it is apparent that $\Gamma_{\lambda}(z)$ is analytic on $\mathbb{C}$ except for the points $z = 0,-1,-2,\dots,\frac{1}{\lambda}, \frac{1}{\lambda}+1,\frac{1}{\lambda}+2, \dots$, where it has simple poles.\\

Thus, by \eqref{21}, we get\\
\begin{equation}\label{22}
\begin{split}
\frac{1}{\Gamma_{\lambda}(z)} &= \lambda^{z}\Gamma\left(\frac{1}{\lambda}\right)z\left(\frac{1}{\lambda}-z\right)\lim_{m \to \infty}e^{\frac{1}{\lambda}\left(1+\frac{1}{2}+\cdots +\frac{1}{m}-\log m\right)}\\
& \quad \times \left[ \lim_{m \to \infty}\prod^{m}_{n=1}e^{-\frac{1}{n\lambda}}\left(1+\frac{z}{n}\right)\left(1+\frac{\frac{1}{\lambda}-z}{n}\right)\right]\\
& = \lambda^{z}z\left(\frac{1}{\lambda}-z\right)\Gamma\left(\frac{1}{\lambda}\right)\lim_{m \to \infty}\left[\prod^{m}_{n=1}m^{-\frac{1}{\lambda}}\left(1+\frac{z}{n}\right)\left( 1+\frac{\frac{1}{\lambda}-z}{n}\right)\right]\\
& = \lambda^{z}z\left(\frac{1}{\lambda}-z\right)\Gamma\left(\frac{1}{\lambda}\right)\lim_{m \to \infty}\left[\prod^{m-1}_{n=1}\left(1+\frac{1}{n}\right)^{-\frac{1}{\lambda}}\right]\\
&\quad \times \left[\prod^{m}_{n=1}\left(1+\frac{z}{n}\right)\left(1+\frac{\frac{1}{\lambda}-z}{n}\right)\right]\\
& =\lambda^{z}z\left(\frac{1}{\lambda}-z\right)\Gamma\left(\frac{1}{\lambda}\right)\prod^{\infty}_{n=1}\left\{\left(1+\frac{1}{n}\right)^{-\frac{1}{\lambda}}\left(1+\frac{z}{n}\right)\left(1+\frac{\frac{1}{\lambda}-z}{n}\right)\right\} .
\end{split}
\end{equation}

Therefore, by \eqref{21} and \eqref{22} we obtain the following theorem.
\begin{thm}\label{Theorem 4} Let $\lambda \in (0,1)$. Then, for all $z \in \mathbb{C}$, with $z \neq  0,-1,-2,\dots,\frac{1}{\lambda}, \\
\frac{1}{\lambda}+1,\frac{1}{\lambda}+2, \dots$, we have
	\begin{equation*}
	\begin{split}
	\Gamma_{\lambda}(z)& = \frac{\lambda^{-z}}{z\left(\frac{1}{\lambda}-z\right)}\frac{1}{\Gamma\left(\frac{1}{\lambda}\right)}e^{-\frac{\gamma}{\lambda}}\prod^{\infty}_{n=1}\left\{ e^{\frac{1}{n\lambda}}\left(1+\frac{z}{n}\right)^{-1}\left( 1+\frac{\frac{1}{\lambda}-z}{n}\right)^{-1}\right\}\nonumber\\
&= \frac{\lambda^{-z}}{z\left(\frac{1}{\lambda}-z\right)}\frac{1}{\Gamma\left(\frac{1}{\lambda}\right)}\prod^{\infty}_{n=1}\left\{\left(1+\frac{1}{n}\right)^{\frac{1}{\lambda}}\left(1+\frac{z}{n}\right)^{-1}\left(1+\frac{\frac{1}{\lambda}-z}{n}\right)^{-1}\right\}.
	\end{split}
	\end{equation*}
\end{thm}

From Theorem 2.4, we note that
\begin{equation}\label{23}
\begin{split}
&\Gamma_{\lambda}(z) =  \frac{\lambda^{-z}}{\Gamma\left(\frac{1}{\lambda}\right)}\frac{1}{z\left(\frac{1}{\lambda}-z\right)}\prod^{\infty}_{n=1}\left\{ \left(\frac{n+1}{n}\right)^{\frac{1}{\lambda}}\left(\frac{n}{n+z}\right)\left(\frac{n}{n+\frac{1}{\lambda}-z}\right)\right\}\\
& =  \frac{\lambda^{-z}}{\Gamma\left(\frac{1}{\lambda}\right)}\frac{1}{z\left(\frac{1}{\lambda}-z\right)}\lim_{m \to \infty}\Bigg\{m^{\frac{1}{\lambda}}\frac{(m-1)!}{(1+z)\cdots (m-1+z)}\\
&\times\frac{(m-1)!}{\left(1+\frac{1}{\lambda}-z\right)\cdots \left(m-1+\frac{1}{\lambda}-z\right)}\Bigg\}
\end{split}
\end{equation}
\begin{equation*}
= \frac{\lambda^{-z}}{\Gamma\left(\frac{1}{\lambda}\right)}\lim_{m \to \infty}\frac{m^{\frac{1}{\lambda}}\left((m-1)!\right)^{2}}{z(1+z)\cdots(m-1+z)\left(\frac{1}{\lambda}-z\right)\left(1+\frac{1}{\lambda}-z\right)\cdots\left(m-1+\frac{1}{\lambda}-z\right)}.
\end{equation*}

Therefore, by \eqref{23}, we obtain the following theorem.\\

\begin{thm}
Let $\lambda \in (0,1)$. Then, for all $z \in \mathbb{C}$, with $z \neq  0,-1,-2,\dots,\frac{1}{\lambda}, \\
\frac{1}{\lambda}+1,\frac{1}{\lambda}+2, \dots$, we have
\begin{equation}\label{Theorem 5}
\Gamma_{\lambda}(z) =  \frac{\lambda^{-z}}{\Gamma(\frac{1}{\lambda})} \lim_{n \rightarrow \infty}
\frac{n^{\frac{1}{\lambda}} ((n-1)!)^2} {z (1+z) \cdots (n-1+z) (\frac{1}{\lambda} -z) (1+\frac{1}{\lambda} -z) \cdots (n-1+\frac{1}{\lambda}-z)}.
\end{equation}
\end{thm}

Let $\lambda \to 1-$ in \eqref{Theorem 5}. Then we have
\begin{equation}\label{24}
\begin{split}
\Gamma_{1}(z) & = \lim_{n \rightarrow \infty} \frac{n ((n-1)!)^2} {z (1+z) \cdots (n-1+z) (1 -z) (2 -z) \cdots (n-1-z) (n-z)} \\
              & = \lim_{n \rightarrow \infty} \frac{((n-1)!)^2} {z (1^2 -z^2) (2^2 -z^2) \cdots ((n-1)^2-z^2)} \\
              & = \frac{1}{z}\lim_{n \rightarrow \infty}  \frac{1} {(1-(\frac{z}{1})^2) (1 -(\frac{z}{2})^2) \cdots (1-(\frac{z}{n-1})^2)} \\	
             & = \frac{1}{z} \prod_{n=1}^{\infty} \bigg( 1-\frac{z^2}{n^2} \bigg)^{-1}. 	
\end{split}
\end{equation}

By Lemma 2.1, and taking $\lambda \to 1-$, we get
\begin{equation}\label{25}
\Gamma_{1}(z) = \frac{\Gamma(z) \Gamma(1-z)}{\Gamma(1)} = \frac{\pi}{\sin \pi z}.
\end{equation}

Therefore, by \eqref{24} and \eqref{25}, we obtain the following corollary.
\begin{cor}
For all $z \in \mathbb{C}-\mathbb{Z}$, we have
\begin{equation*}\label{cor 6}
\frac{\pi z}{\sin \pi z} = \prod_{n=1}^{\infty} \bigg( 1-\frac{z^2}{n^2} \bigg)^{-1}.
\end{equation*}
\end{cor}

In light of the equation \eqref{03}, it is natural to define the degenerate beta functions by
\begin{equation}\label{25-1}
B_{\lambda}(\alpha,\beta)=\frac{\Gamma_{\lambda}(\alpha)\Gamma_{\lambda}(\beta)}{\Gamma_{\lambda}(\alpha+\beta)},
\end{equation}
where $\alpha, \beta \neq 0,-1,-2,...\frac{1}{\lambda},\frac{1}{\lambda}+1,\frac{1}{\lambda}+2,...$.\\

From Lemma 2.1, Theorem 2.3  and the first identity in Theorem 2.4, it is immediate to derive the following expressions for the degenerate beta function. \\
\begin{thm}
Let $\lambda \in (0,1)$. Then, for any $\alpha, \beta  \neq 0,-1,-2,...\frac{1}{\lambda},\frac{1}{\lambda}+1,\frac{1}{\lambda}+2,...$, we have the following expressions.
\begin{equation*}\begin{split}
B_{\lambda}(\alpha,\beta)
& = \frac{1}{ \Gamma(\frac{1}{\lambda})} B(\alpha, \beta) \frac{\Gamma(\frac{1}{\lambda} - \alpha)  \Gamma(\frac{1}{\lambda} - \beta)} {\Gamma(\frac{1}{\lambda} -\alpha -\beta) }\\
&=\frac{e^{-\frac{\gamma}{\lambda}}(\alpha+\beta)(\frac{1}{\lambda}-\alpha-\beta)}{\Gamma(\frac{1}{\lambda})\alpha\beta(\frac{1}{\lambda}-\alpha)(\frac{1}{\lambda}-\beta)}\\
&\times\quad\prod_{n=1}^{\infty}\frac{e^{\frac{1}{n\lambda}}(1+\frac{\alpha+\beta}{n})(1+\frac{\frac{1}{\lambda}-\alpha-\beta}{n})}{(1+\frac{\alpha}{n})(1+\frac{\frac{1}{\lambda}-\alpha}{n})(1+\frac{\beta}{n})(1+\frac{\frac{1}{\lambda}-\beta}{n})}.
\end{split}\end{equation*}
In particular, if $\alpha=m, \beta=n$ are positive integers, then we have the following expression.
\begin{equation*}
B_{\lambda}(m,n)=  \frac{(1)_{m+n+1,\lambda}}{(1)_{m+1,\lambda} (1)_{n+1,\lambda} } B(m, n).
\end{equation*}
\end{thm}

\vspace{1mm}

\section{\bf Further Remark}
\vspace{0.05 in}

As is traditional, let $\int_{\infty}^{(0+)}$ (respectively, $\int_{-\infty}^{(0+)})$ denote the integration of path that starts at infinity (respectively, at negative infinity), encircles the origin counter-clockwise direction and returns to the starting point. Assume that $0<\lambda<1, \, 0<\textnormal{Re}(s)<\frac{1}{\lambda}$.\\
\indent As we noted in \eqref{16}, we have
\begin{equation}\label{26}
\lambda^s\Gamma_{\lambda}(s)=\int_{0}^{\infty}t^{s-1}(1+t)^{-\frac{1}{\lambda}}dt.
\end{equation}
We now consider the integral given by
\begin{equation}\label{27}
\int_{C}\frac{z^{s-1}}{(1-z)^{\frac{1}{\lambda}}}dz.
\end{equation}
Here $C$ is the integration of path that moves along the negative real axis from $-R$ to $-\delta$, along the circle of radius $\delta \, (<1)$ in counter-clockwise direction and then along the negative real axis from $-\delta$ to $-R$. Then the integral in \eqref{27} is equal to
\begin{equation}\label{28}
\int_{R}^{\delta}\frac{t^{s-1}e^{-i \pi s}}{(1+t)^{\frac{1}{\lambda}}}dt
+\int_{-\pi}^{\pi}\frac{iR^se^{is\theta}}{(1-Re^{i \theta})^{\frac{1}{\lambda}}} d\theta
+\int_{\delta}^{R}\frac{t^{s-1}e^{i \pi s}}{(1+t)^{\frac{1}{\lambda}}}dt.
\end{equation}
\indent Let $R \rightarrow \infty$, and let $ \delta \rightarrow 0$. Then the second integral in \eqref{27} tends to zero and hence we get
\begin{equation}\label{29}
\lambda^s\Gamma_{\lambda}(s)=\int_{0}^{\infty}t^{s-1}(1+t)^{-\frac{1}{\lambda}}dt
=\frac{1}{2 i \sin \pi s}\int_{-\infty}^{+0} \frac{z^{s-1}}{(1-z)^{\frac{1}{\lambda}}} dz.
\end{equation}
\indent We note here that the rightmost integral in \eqref{29} vanishes at $s=1,2,3,\dots$ and it does not vanish at  $s=0,-1,-2,-3,\dots$. Thus this gives an analytic continuation $\Gamma_{\lambda}(s)$ to a meromorphic function on $\mathbb{C}$, with simple poles at $s=0,-1,-2,-3,\dots$.
Further, in view of the symmetry of  $\lambda^{-s}\Gamma_{\lambda}(s)$ under $s \rightarrow \frac{1}{\lambda}-s$, it also has simple poles at $\frac{1}{\lambda},\frac{1}{\lambda}+1,\frac{1}{\lambda}+2, \dots$. As we can see from Theorem 2.4, $s=0,-1,-2,-3,\dots, \frac{1}{\lambda},\frac{1}{\lambda}+1,\frac{1}{\lambda}+2, \dots$ are all the poles of $\Gamma_{\lambda}(s)$.\\
\indent Finally, by replacing $z$ by $-z$ in \eqref{29} we also have the expression given by
\begin{equation}
\lambda^s\Gamma_{\lambda}(s)
=\frac{i}{2  \sin \pi s}\int_{\infty}^{+0} \frac{(-z)^{s-1}}{(1+z)^{\frac{1}{\lambda}}} dz.
\end{equation}

\vspace{.5cm}

\section{\bf Conclusion}

\medskip

The degenerate Bernoulli and Euler polynomials were introduced by Carlitz as degenerate versions of ordinary Bernoulli and Euler polynomials, respectively. In recent years, studying degenerate versions of some special polynomials and numbers regained interests of many mathematicians and truned out to be very useful and fruitful. This idea of investigating various degenerate versions is not just limited to polynomials but can be extended to transcendental functions. Indeed, the degenerate gamma functions were introduced in [6] as a degenerate version of the usual gamma function (see also [7]). \\
\indent  In this paper, we derived several basic properties about the degenerate gamma functions. Firstly, we noted that, by making use of \eqref{16}, the degenerate gamma function $\Gamma_{\lambda}(s)$ initially defined for $0 < \textnormal{Re}(s) <\frac{1}{\lambda}$ can be analytically continued to a meromorphic function on $\mathbb{C}$, except for simple poles at $s=0,-1,-2,...\frac{1}{\lambda},\frac{1}{\lambda}+1,\frac{1}{\lambda}+2,...$.  Secondly, we derived the difference formula in Theorem 2.2, determined their values at positive integers in Theorem 2.3, and several expressions by making use of the Weierstrass and Euler formulas for the usual gamma function in Theorems 2.4 and 2.5. Thirdly, we defined the degenerate beta functions and obtained some expressions for them. Lastly, we obtained an integral representation for the degenerate gamma functions as an integral along a Hankel contour.








\end{document}